\documentclass[a4paper, 11pt]{article}
\usepackage{amsmath} \usepackage{amsfonts} \usepackage{amssymb}\usepackage{latexsym}
\usepackage[english]{babel}
\usepackage{amscd}
\usepackage[dvips,final]{graphics}
\usepackage{locengli,math}
\usepackage{graphicx,pst-plot,psfig,pstricks,pst-node,
pst-text,pst-tree,pst-char,pst-coil,Addaletree}
\newcommand{\N}{\hat{\mathbb{N}}}

\newcommand{\Na}{\mathbb{N}}

\def\TT{T\!\!\!\! T}
\def\NN{N\!\!\!\! N}

\def\u{\underline }
\newcommand{\Ng}{\hat{\NN}}

\newcommand{\Co}{ \textsf{Corl.}}

\begin{document}
\begin{center}
\textbf{\LARGE{\textsf{A $LL$-Lattice reformulation of arithmetree over planar rooted trees. Part II}}}
\footnote{Supported by the European Commission HPRN $\sim$ CTN2002 $\sim$ 00279, RTN QP-Applications.
{\it{2000 Mathematics Subject Classification:   05C05; 06A07; 11A99; 06A07.}}
{\it{Key words and phrases: $LL$-lattice, dendriform trialgebra over one generator, planar rooted trees, arithmetree, involutive $\mathcal{P}$-Hopf algebras.}} 
}
\vskip1cm
\parbox[t]{14cm}{\large{
Philippe {\sc Leroux}}\\
\vskip4mm
{\footnotesize
\baselineskip=5mm
Institut f\"ur Mathematik und Informatik,\\
Ernst-Moritz-Arndt-Universit\"at, Jahnstra$\beta$e 15a, 17487 Greifswald, Germany,\\ leroux@uni-greifswald.de}}
\end{center}
\vskip1cm
\baselineskip=5mm
\noindent
\begin{center}
\begin{large}
\textbf{ 30/06/04}
\end{large}
\end{center}
{\bf Abstract:} 
We continue our reformulation of free dendriform algebras, dealing this time with the free dendriform trialgebra generated by $\treeA$ over planar rooted trees. We propose a `deformation' of a vectorial coding used in \cite{Lerden}, giving a $LL$-lattice on rooted planar trees according to the terminology of A. Blass and B. E. Sagan. The
three main operations on trees become explicit, giving thus a complementary approach to a very recent work of P. Palacios and M. Ronco. Our parenthesis framework allows a more tractable reformulation to explore the properties of the underlying lattice describing operations and simply a proof of a fundamental theorem related  
to arithmetics over trees, the so-called arithmetree. Arithmetree is then viewed as a noncommutative extention of $(\Na,+,\times)$, the integers being played by the corollas.
\section{Introduction}
In the sequel, $K$ is a null characteristic field and $\Na$ is the semiring of integers. If $S$ is a finite set, then $\card(S)$ denotes its cardinal and $KS$, the $K$-vector space spanned by $S$.
Rooted planar trees, often called trees for short, are known to be in bijection with all possible parentheses constructed over $x_1, \ldots, x_n$ and modeling at least binary operations. In the sequel, by complete expression, we mean a monomial
of $\bra \ x_1, \ldots, x_n, ( ,)\ \ket$ --the free associative semigroup generated by $x_1, \ldots, x_n, ($ and  $)$-- in one-to-one correspondence with a rooted planar tree, \textit{i.e.,} every $($ is closed by a unique $)$.
In \cite{Lerden}, we proposed a reformulation of the dendriform dialgebra on the generator $\treeA$ over rooted planar binary trees \textit{via} a parenthesis framework.  Complete expressions of $\bra x_1, \ldots, x_n, (, ) \ket$ were canonically associated with rooted planar binary trees obtaining thus an injection map $ Exp:   Y_n \xrightarrow{} \bra \ x_1, \ldots, x_{n+1}, ( ,) \ \ket$, where $Y_n$ is the set of rooted planar binary trees with $n$ internal vertices.
Parentheses of a complete expression
were coded into a unique vector  
of $\Na^{n}$ obtained as follows. 
Encode the parentheses of $Exp(\tau)$ of the binary tree $\tau$ in a vector  $\vec{v} :=(v_1, v_2, \ldots, v_{n})$ of $\Na^{n}$ by declaring that for all $1 \leq i \leq n$,
$v_i :=i$ if and only if there exists a left parenthesis at the left hand side of $x_i$, \textit{i.e.}, $\ldots(^p x_i \ldots$, with $p>0$, occurs in  $Exp(\tau)$. Otherwise, there exists a unique most right parenthesis at the right hand side of $x_i$ which closes a unique left parenthesis say open at $x_j$. In this case, $v_i :=j$. We then obtain an injective map:
$ \vec{name}:   Y_n \xrightarrow{} \Na^{n},$
which map any tree $\tau$ into a vector, $\vec{name}(\tau)$, also denoted by $\vec{\tau}$ for short, called the name of $\tau$.
In the sequel, $\vec{name}(Y_n)$ will be denoted by $ \N^{n}$.
Still in the case of planar binary trees, it was shown \cite{LRbruhat, Lerden} that once the sets $Y_n$ (or equivalently $\N^n$) where equipped with the Tamari partial order, resp. with the trivial partial order,
operations defining the free dendriform dialgebra over the generator $\treeA$ turned out to be explicit.
We keep this method to deal with planar rooted trees. Contrary to the binary case, the difficulty is now
to keep track of every parenthesis of the form $)$. This will lead to a `deformation' of the sets $\N^{n}$ in Subsections 2.1 and 2.2. These sets
becomes lattices when equipped with a very natural partial order called trivial partial order.  In Subsection 2.2, we study how to construct the meet and the joint
of two elements. We proved that the dendriform involution is a lattice anti-automorphism and propose a way to compute names of trees subject to involution. Using results of B.E. Sagan and A. Blass  \cite{Sagan, BlassSagan}, we compute the M\"obius function of these lattices and proved they are of type $LL$ (like for the Tamari lattices $Y_n$ or $\N^n$ associated with planar binary trees). Characteristic polynomials of these lattices are computed. Furthemore, we rediscover as a consequence of our vectorial coding, the definition of 
a partial order introduced very recently by P. Palacios and M. Ronco \cite{PalRon}. 
In Subsection 2.3, we propose another framework, based on this parenthesis point of view, to explicit the operations defining the dendriform trialgebra over the generator $\treeA$. In Subsection 2.4, we propose a lattice formulation of arithmetree over planar rooted trees by notably simplifying a fundamental proof due to J.-L. Loday \cite{Lodayarithm}. We explicit the coproduct of an involutive Hopf algebra associated with planar rooted trees and show that
trees endowed with their arithmetree can be viewed as a noncommutative version of our usual arithmetics over integers. More precisely, we construct a Hopf algebra over $(K\N, +)$ and establish
an isomorphism of associative algebras between $(K\N, +, \times)$
and the associative algebra generated by the corollas. On $(K\N, +)$, the generator is known to be $[1]$ which is mapped \textit{via} this isomorphism
to the corolla $\treeA$. However, the associative operation $+$ on integers has to be replaced by three operations on planar rooted trees compatible with the action of the neutral element denoted by $[0]$. This replacement will generate a modification of the structure of the Hopf algebra on integers to give an involutive one computed in this paper. We conclude by enumerating planar trees, noncommutative generalisation of our usual integers, invariant under the dendriform involution. We also find
two new interpretations of the super Catalan numbers.
\section{Rooted planar trees}
Denote by $T_n$
the set of rooted planar trees with $n+1$ leaves, \textit{i.e,} one root and each internal vertex with at least two leaves. Consider them up to isotopies. In small dimension, we obtain:
$$ T_0:= \{ (0):=\treeO \}, \  T_1:=\{ (1):=\treeA \}, \ T_2:=\{ \treeAB, \ \treeM, \  \treeBA \},$$ $$T_3:=\{ \treeABC, \treeMB, \treeBAC, \treeAMC,  \treeACA, \treeCAB,  \treeAM, \treeMA,  \treeCM, \treeCBA, \treeCor \}.$$
The cardinal of the $T_n, n>0,$ are the super Catalan numbers or Schr\"oder numbers and are denoted by $C_n$, \textit{i.e.,} $C_0=1, C_1=1, C_2=3, C_3=11, C_4=45, \ldots$. The grafting operation is still denoted by $\vee$. Every tree $t$ can be uniquely written as 
$t_1 \vee t_2 \vee \ldots \vee t_n$ where the $t_i$ are also trees. Pictorially, the roots of the $t_i$ are glued together, forming a unique root, the root of $t$. Example: $\treeMA:= \treeO \vee \treeO \vee \treeA$. The elements $\Co_{[n]}:=\treeO \vee  \treeO \vee \ldots \vee \treeO$, $n+1$ times are called corollas. There exists an involution called the dendriform involution defined inductively by $t^\dagger:=t_n^\dagger \vee t_{n-1}^\dagger \vee \ldots \vee t_1^\dagger$, if $t:=t_1 \vee t_2 \vee \ldots \vee t_n$ is a planar tree.
This  dendriform involution will play an important r\^ole in the sequel of this paper.
In \cite{LodayRonco}, M. Ronco and J.-L. Loday introduced dendriform trialgebras which are $K-$vector spaces
$T$ equipped with three binary operations:
$\prec, \ \succ, \ \bullet: T^{\otimes 2} \xrightarrow{} T$, satisfying the following relations for all $x,y,z \in T$:
$$(x \prec y )\prec z = x \prec(y \star z), \ 
(x \succ y )\prec z = x \succ(y \prec z), \
(x \star y )\succ z = x \succ(y \succ z), $$
$$(x \succ y )\bullet z = x \succ(y \bullet z), \ 
(x  \prec y )\bullet z = x \bullet(y \succ z), \ 
(x \bullet y )\prec z = x \bullet(y \prec z), \
(x \bullet y ) \bullet z = x \bullet(y \bullet z),$$
where by definition $x \star y :=x  \prec y +x \succ y +x \bullet y$, for all $x,y \in T$, turns out to be associative. This defines a regular, binary and quadratic operad whose Poincar\'e series starts with $1,3, 11 \ldots$ like the super Catalan series. They showed that the augmented free 
dendriform trialgebra on one generator $x$ is isomorphic to 
$KT_{\infty}^*:= \oplus_{n >0} KT_n$, the generator $x$ being mapped to the generator $\treeA$. The tree $\treeO$ is the unit for the operation $\star$, \textit{i.e,} $ t \star \treeO := t =:   \treeO \star t $ and the operations $\prec, \ \succ, \ \bullet$ are given on trees  inductively by the following formulas:  for any trees $t:=t_1 \vee \ldots \vee t_n$ and $z:=z_1 \vee \ldots \vee z_m $,
\begin{enumerate}
\item {$t\succ z:= (t \star z_1) \vee z_2 \vee \ldots \vee z_m$,}
\item {$t\prec z:= t_1 \vee t_2 \vee \ldots \vee (t_n \star z)$,}
\item {$t \bullet z:= t_1 \vee t_2 \vee \ldots \vee (t_n \star z_1) \vee z_2 \vee \ldots \vee z_m$.}
\end{enumerate} 
The aim of this section is to code rooted planar trees, to propose a natural partial order over them, generalising the Tamari one, and to give a complementary point of view to a very recent work of P. Palacios and M. Ronco \cite{PalRon}. 
\subsection{Deformation of $K\Na^n$}
We fix here some useful notation.
Fix $n>0$. Consider the set $\Na[h^{-1}]^{n}$ consisting of vectors $\vec{v}:=(v_1[h^{-1}], v_2[h^{-1}],\ldots,v_{n}[h^{-1}] )$, where the $v_i[h^{-1}] \in \Na[h^{-1}]$, \textit{i.e.,} are polynomials in $h^{-1}$ with coefficients in $\Na$. The degree of a polynomial from $\Na[h^{-1}]^{n}$
is the absolute value of the power of its lower monomial. 
The so-called \textit{trivial partial order} on $\Na[h^{-1}]$ is induced by the following rule.

\noindent
\textbf{Rule:} Let $P:= a_0 + a_1h^{-1} + \ldots + a_nh^{-n}$ and $Q:=b_0 + b_1h^{-1} + \ldots + b_mh^{-m} \in \Na[h^{-1}]$. Then $P>Q$ if and only if there exists a $0 \leq k \leq \max(n,m)=\max(\deg(P), \deg(Q)$ such that for all $0 \leq i<k$, $a_i=b_i$ and $a_k>b_k$ (the case $k=0$ meaning that $a_0>b_0$).

\noindent
This induces the so-called \textit{trivial partial order} on  $\Na[h^{-1}]^n$, for all $n>0$, by declaring that $\vec{v}:=(v_1[h^{-1}], v_2[h^{-1}],\ldots,v_{n}[h^{-1}] ) < \vec{w}:=(w_1[h^{-1}], w_2[h^{-1}],\ldots,w_{n}[h^{-1}] )$ if and only if
for all $1 \leq i \leq n$, $v_i[h^{-1}] \leq w_i[h^{-1}]$, with at least one strict inequality. 
The introduction of the set $\Na[h^{-1}]^{n}$ is motivated by \cite{Lerden}. Recall, see the introduction, that every rooted planar binary tree was coded into a vector with integer coordinates and that only most right parentheses were sufficient to code $\ldots x_k)^i \ldots$. For planar trees, we have to take into account all the closing parentheses. Hence, instead of natural numbers and their usual order, we consider an analogue of real numbers in some sense, usually written as $a_0 + a_1 10^{-1} +
a_2   10^{-2} + \ldots  $ and their usual order defined as above where $h$ has to be replaced by 10 in this particular case. 
\subsection{A lattice on rooted planar trees}
We now propose a natural way to code parentheses by vectors. Fix $n>0$. With any $t \in T_n$, we will associate a unique vector of $\Na[h^{-1}]^{(n+1)}$. 
Any planar tree $t$ defines a unique complete expression $Exp(t)$ in $\bra x_1, \ldots, x_n,x_{n+1}, (,) \ket$ and hence a unique vector $\vec{v}:=(v_1[h^{-1}], \ldots,v_{n+1}[h^{-1}]) \in \Na[h^{-1}]^{(n+1)}$ defined as follows.  Fix $1\leq i \leq n+1$. The coordinate $v_i[h^{-1}]:=i$, if there exists $(^p$, $p>0$, at the left hand side of $x_i$, \textit{i.e.,} locally the monomial $Exp(t)$ has the form $\ldots (^p x_i  \ldots $.  The coordinate $v_i[h^{-1}]:=i-1 + ih^{-1}$ if locally the monomial $Exp(t)$ has the form $\ldots )^p x_i (^q  \ldots $ with $p,q \in \Na$.  The coordinate $v_i[h^{-1}]:=p_1+h^{-p_1} + \ldots + h^{-p_k}$, where $1 \leq p_1 \leq \ldots \leq p_k < n,$ if locally the monomial $Exp(t)$ has the form $\ldots  x_i)^k  \ldots $ with $k>0$ and closing parentheses open at $x_{p_1}, \ldots ,x_{p_{k-1}}$ and $x_{p_{k}}$. Observe that the most right parenthesis in $\ldots  x_i)^k  \ldots $ closes one $($ in $x_{p_1}$, fixing so the null degree part of $v_i[h^{-1}]$. The unique vector representing a tree is still called its name and we identify $T_n$
with $\N[h^{-1}]^{n+1}$, the set of vectors naming trees from $T_n$.
The trivial partial order defined on $\Na[h^{-1}]^{n+1}$ induces a trivial partial order on $\N[h^{-1}]^{n+1}$ and thus on $T_n$.  Here are examples in small dimensions. By convention $\treeO:= (0)$. The tree $\treeA$ is equal to $(1, 1+h^{-1})$. For esthetic reasons, if $t,t'$ are trees, then $t \rightarrow t'$ will be equivalent to $t < t'$.
The lattice $(T_2, <)$ or $(\N[h^{-1}]^2,<)$ is of the form,
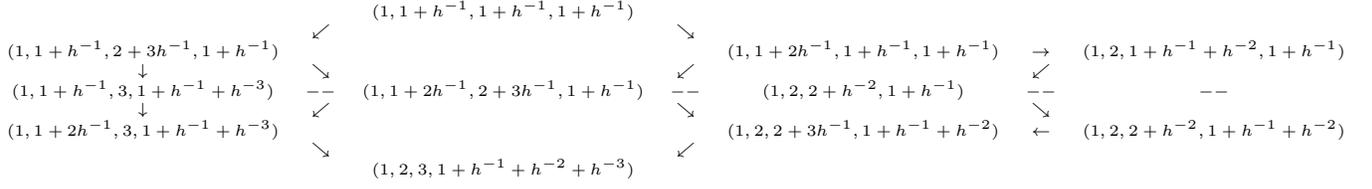
$$ \treeAB \longrightarrow \treeM \longrightarrow \treeBA,$$
\begin{footnotesize}
$$ (1,1+h^{-1},1+h^{-1}) \longrightarrow (1,1+2h^{-1},1+h^{-1}) \longrightarrow (1,2,1+h^{-1}+h^{-2}).$$
\end{footnotesize}
and $(T_3, <)$ or $(\N[h^{-1}]^3, <)$ is as follows.
\begin{center}
$
\begin{array}{ccccccccccc}
 & & & \treeABC & & & & & & & \\
 & & \swarrow & & \searrow & & & & & & \\
 & \treeAM &  & &  & \treeMB & \longrightarrow & \treeBAC & & & \\
 & \downarrow & \searrow & &\swarrow  & &  & &\searrow & & \\
 -- & \treeACA &--& \treeCor  &-- & -- &-- & --&-- &\treeAMC & --\\
& \downarrow & \swarrow & & \searrow & & & &\swarrow & & \\
& \treeMA &  & &  & \treeCM &\longleftarrow  & \treeCAB & & & \\
& & \searrow & & \swarrow & & & & & & \\
& & & \treeCBA & & & & & & & \\
\end{array}
$
\end{center}
where the dash line is the symmetry axis determined by the dendriform involution.
\begin{center}
\begin{tiny}
$
\begin{array}{ccccccccc}
  & & (1,1+h^{-1},1+h^{-1},1+h^{-1}) & & & & & &  \\
  & \swarrow & & \searrow & & & & &  \\
 (1,1+h^{-1}, 2+3h^{-1},1+h^{-1}) &  & &  & (1,1+2h^{-1},1+h^{-1},1+h^{-1}) & \rightarrow & (1,2,1+h^{-1}+h^{-2},1+h^{-1}) & &  \\
  \downarrow & \searrow & &\swarrow  & & \swarrow & & &  \\
 (1,1+h^{-1},3,1+h^{-1}+h^{-3}) &--& (1,1+2h^{-1},2+3h^{-1},1+h^{-1})  &-- &(1,2,2+h^{-2},1+h^{-1})  &-- & --&  &  \\
 \downarrow & \swarrow & & \searrow & &\searrow & & &  \\
 (1,1+2h^{-1},3,1+h^{-1}+h^{-3}) &  & &  & (1,2,2+3h^{-1},1+h^{-1}+h^{-2}) &\leftarrow  &(1,2,2+h^{-2},1+h^{-1}+h^{-2}) & &  \\
 & \searrow & & \swarrow & & & & &  \\
 & & (1,2,3,1+h^{-1}+h^{-2}+h^{-3}) & & & & & &  \\
\end{array}
$
\end{tiny}
\end{center}
\begin{prop}(Involution)
Fix $\vec{v}[h^{-1}] \in \N[h^{-1}]^{n}$, $n>0$. Then, the vector $\vec{v}[h^{-1}]^\dagger$ is obtained as follows.
\begin{enumerate}
\item {If $v[h^{-1}]_i :=i$, then  $v[h^{-1}]^\dagger _{n+1-i} :=(n+1-i_k)+h^{-(n+1-i_k)}+ \ldots h^{-(n+1-i_0)}$, where the $i_0 < i_1< \ldots <i_k$ are the positions where $h^{-i}$ appears in $\vec{v}[h^{-1}]$.}
\item {If $v[h^{-1}]_i :=i-1+ih^{-1}$, then  $v[h^{-1}]^\dagger _{n+1-i} :=n-i+(n+1-i)h^{-1}$.}
\item {If $v[h^{-1}]_i :=i_0+h^{-i_0}+ \ldots +h^{-i_k}$, then  $v[h^{-1}]^\dagger _{n+1-i} :=n+1-i$.}
\end{enumerate} 
\end{prop}
\Proof
Observe that the dendriform involution of a complete expression associated with a tree is obtained by reading it from left to right, the closing parentheses becoming open ones and conversely.
\eproof
\begin{prop}(Lattice anti-automorphism)
\label{antipla}
Fix $n>0$ and $\vec{v}[h^{-1}], \vec{w}[h^{-1}] \in \N[h^{-1}]^{n}$. Then,
$$\vec{v}[h^{-1}] < \vec{w}[h^{-1}] \Leftrightarrow \vec{w}[h^{-1}]^\dagger < \vec{v}[h^{-1}]^\dagger.$$
\end{prop}
\Proof
Fix $n>0$ and $\vec{v}[h^{-1}], \vec{w}[h^{-1}] \in \N[h^{-1}]^{n}$, such that
the left hand side of the previous inequality holds.
We prove Proposition \ref{antipla} by checking every case, which is straightforward except the case where $v[h^{-1}]_i:=i$ and
$w[h^{-1}]_i:=i$, for some $i>1$. We will show in this case that $v^\dagger[h^{-1}]_{N+1-i} \geq w^\dagger[h^{-1}]_{N+1-i}$. For that, we have to focus on the first parenthesis, standing at $i_0>i$, closing one open in $i$. Suppose the existence of a $)$, standing at $i_1$, between $i$ and $i_0$ in the expression associated with $\vec{w}[h^{-1}]$. This will imply the existence of a $)$, standing at $i_1$, in the expression associated with $\vec{v}[h^{-1}]$ and closing one $($ open in $i_2$, with $i < i_2 <i_1$, by hypothesis. Observe then that $w[h^{-1}]_{i_1}:=i_1'+ \ldots $ where $i_1' \leq i$ is the $($ standing at $i_1'$ and closed by the most external $)$ in $i_1$ whereas $v[h^{-1}]_{i_1}:=i_1'' + \ldots$, where $i_1'' \geq i$. Hence, the first parenthesis in the expression of  $\vec{w}[h^{-1}]$ closing one $($ open in $i$ has to be placed after or at the same position that the first parenthesis in the expression of  $\vec{v}[h^{-1}]$ closing one $($ open in $i$. By repeating this remark, in case of equality in the positions of parentheses, we get, $v^\dagger[h^{-1}]_{N+1-i} \geq w^\dagger[h^{-1}]_{N+1-i}$.
By checking every case, one obtains $v^\dagger[h^{-1}]_{N+1-i} \geq w^\dagger[h^{-1}]_{N+1-i}$, for all $1 \leq i \leq n$. Hence, $\vec{v}^\dagger[h^{-1}] \geq \vec{w}^\dagger[h^{-1}]$. But the case $\vec{v}^\dagger[h^{-1}] = \vec{w}^\dagger[h^{-1}]$ is impossible since the dendriform involution is involutive.
\eproof
\begin{prop}(Meet-joint)
\label{meetjoint}
Fix $n>0$ and $\vec{v}[h^{-1}], \vec{w}[h^{-1}] \in \N[h^{-1}]^{n}$. Then,
their least upper bound, or joint denoted by $\vec{lup}[h^{-1}]:= \bigvee(\vec{v}[h^{-1}], \vec{w}[h^{-1}])$ is obtained as follows.
\begin{enumerate}
\item {$lup[h^{-1}]_i:=i$ if one of the $v[h^{-1}]_i$ or $w[h^{-1}]_i$ is equal to $i$.}
\item {$lup[h^{-1}]_i:=i-1+ih^{-1}$ if one of the $v[h^{-1}]_i$ or $w[h^{-1}]_i$ is equal to $i-1+ih^{-1}$, the other being least or equal to it.}
\item {Suppose $v[h^{-1}]_i:= i_0 + h^{-i_0}+ \ldots + h^{-i_k}$ and 
$v[h^{-1}]_i:= j_0 + h^{-j_0}+ \ldots + h^{-j_p}$, where both $i_0$ and $j_0$ are less than $i$, and, $p$ and $k$ are integers. Then,
$lup[h^{-1}]_i:=a_0+h^{-a_0}+ \sum_{r \in I} h^{-r}$ where  $a_0:=\max(i_0,j_0)$ and $I$ is the set of labels between $a_0$ and $i$ of
parentheses $($, --\textit{i.e.}, locally in the expression associated with $\vec{lup}[h^{-1}]$, we have $...(X_r...$-- which have not been closed by a $)$ yet.}
\end{enumerate}
Their greatest lower bound or meet, denoted by $\vec{glw}[h^{-1}]:= \bigwedge(\vec{v}[h^{-1}], \vec{w}[h^{-1}])$ is obtained as follows.
\begin{enumerate}
\item {$glw[h^{-1}]_i:= min(v[h^{-1}]_i, w[h^{-1}]_i)$, except if one of them is equal or less than $i-1 + ih^{-1}$.}
\item{If say, $v[h^{-1}]_i \geq i-1 + ih^{-1}$ and $w[h^{-1}]_i:= j_0 + h^{-j_0}+ \ldots + h^{-j_p}$, then $glw[h^{-1}]_i:=  a_0 + h^{-a_0}+ \sum_{r \in I} h^{-r}$, where $a_0:=\max \{k, v[h^{-1}]_k=k=w[h^{-1}]_k \}$ and $a_0 \leq j_0$ and where $I$ is the set between $a_0$ and $i$ of
parentheses $($ which have not been closed by a $)$ yet.}
\item{If $v[h^{-1}]_i:=  i_0 + h^{-i_0}+ \ldots + h^{-i_k}$ and $w[h^{-1}]_i:= j_0 + h^{-j_0}+ \ldots + h^{-j_p}$, then $glw[h^{-1}]_i:=  a_0 + h^{-a_0}+ \sum_{r \in I} h{-r}$, where $a_0:=\max \{k, v[h^{-1}]_k=k=w[h^{-1}]_k \}$ and $a_0 \leq min(i_0, j_0)$ and where $I$ is the set between $a_0$ and $i$ of
parentheses $($ not closed by a $)$ yet.}
\end{enumerate} 
Moreover, the meet and the joint are related as follows $(\bigvee(\vec{v}[h^{-1}],  \vec{w}[h^{-1}]))^\dagger = \bigwedge(\vec{v}[h^{-1}]^\dagger, \vec{w}[h^{-1}]^\dagger)$.
\end{prop}
\Proof
Keep notation of this Proposition. 
Conditions 1 and 2 are obvious. Condition 3, means that the most external parenthesis in the expression of $\vec{v}[h^{-1}]$ (resp. $\vec{w }[h^{-1}]$) closes a $($ open in $i_0$ (resp. $j_0$). Hence, the least upper boud has to have in coordinate $i$, its most external parenthesis $)$, closing one open in $a_0:=\max(i_0,j_0)$. To avoid contradiction, we have to close every $($ between $a_0$ and $i$ in the expression associated with $\vec{lup}[h^{-1}]$, hence Item 3. The proof is the same for the meet.
The last claim holds by considering the dual lattice obtained under the action of the dendriform involution which is a lattice anti-automorphism. 
\eproof

\noindent
To state Theorem \ref{Sagan}, introduce the set $\mathcal{A}(L)$ of a lattice $L$ with minimal element, to denote the set of all atoms of $L$, ---those elements such that there is no other one between them and the minimum---. Such a set is call \textit{independent} if for all $B \subsetneq \mathcal{A}(L)$, $\bigvee B < \bigvee \mathcal{A}(L)$, where $\bigvee$ stands for the least upper bound operation. The following result holds. 
\begin{theo}[B. E. Sagan \cite{Sagan}]
\label{Sagan}
Let $L$ be a finite lattice such that $\mathcal{A}(L)$ is independent. Then, the M\"obius function $M$ of $L$ is
$ M(x)= (-1)^{\card B}$ if $x:=\bigvee B$, for some $B \subseteq \mathcal{A}(L)$, and $ M(x)=0$ otherwise.
\end{theo}
\begin{coro}
Denote by $M_h$, the M\"obius function of the lattice $(\N[h^{-1}]^n, <)$. Fix $\vec{v}[h^{-1}] \in \N[h^{-1}]^{n}$. Then, $M_h(\vec{v}[h^{-1}])=(-1)^t$, if every coordinate $v[h^{-1}]_i$ is equal to either $1+h^{-1}$ or $i-1+ih^{-1}$. In this case $t$ is number of $v_i \not=1+h^{-1}$. Otherwise, $M_h(\vec{v}[h^{-1}])=0$.
\end{coro}
\Proof
Observe that the set of atoms $\mathcal{A}(\N[h^{-1}]^n, <)$ of this lattice is independent since the $n-1$ atoms are of the form $a_1:=(1, 1+2h^{-1},1+h^{-1}, \ldots, 1+h^{-1}), a_2:=(1, 1+h^{-1}, 2+3h^{-1},1+h^{-1}, \ldots, 1+h^{-1}), \ldots, a_{n-1}:=(1,1+h^{-1}, \ldots, 1+h^{-1}, (n-1)+nh^{-1},1+h^{-1})$. We have $\bigvee(a_1, \ldots,a_{n-1}):= \Co_n$ and for all $B \subsetneq \mathcal{A}(\N[h^{-1}]^n, <)$, $\bigvee B < \Co_n$. Apply Theorem \ref{Sagan} to conclude.
\eproof
\begin{prop}
Fix $n>0$ and  $\vec{v}[h^{-1}] \in \N[h^{-1}]^{(n+1)}$.  Consider the following transformations. 
Replace each $v_i:= i-1 + ih^{-1}$ by $v_i:= i$ and each $v_i:= i_0 + h^{- i_0}
+ h^{- i_1}+ \dots + h^{- i_k}$ by $v_i:= i_0$ and forget the last coordinate. This defines a surjection $P: \N[h^{-1}]^{(n+1)} \rightarrow \N^{n}$.
If $\vec{v}[h^{-1}], \vec{w}[h^{-1}] \in \N[h^{-1}]^{(n+1)}$, then $\vec{v}[h^{-1}] < \vec{w}[h^{-1}] \Rightarrow P(\vec{v}[h^{-1}]) \leq P(\vec{w}[h^{-1}])$, the inequality being strict if both vectors named rooted binary trees. Moreover, if $M_h(\vec{v}[h^{-1}]) \not=0$, then $M(P(\vec{v}[h^{-1}])) \not=0$, where $M$ is the M\"obius function on $\N^{n}$.
\end{prop}
\Proof
Fix $n>0$, and $\vec{v}[h^{-1}] \in \N[h^{-1}]^n$. Under the transformation $P$, the last coordinate of $\vec{v}[h^{-1}]$ will give always 1 and for all $i$, $v[h^{-1}]_i< w[h^{-1}]_i$ implies that the restriction of $P$ to the coordinate $i$, still denoted by $P$, will give $P(v[h^{-1}]_i) \leq P(w[h^{-1}]_i)$ in $\Na^{n}$. These vectors give names of rooted planar binary trees. Indeed, $P$ maps
$\ldots X_i \ldots $, belonging to
a complete sub-expression within that associated with $v[h^{-1}]$ or $w[h^{-1}]$, to $\ldots (X_i \ldots $, giving thus a sub-expression. But, there exists a unique way to complete it, by placing $)$ between the last $X$ of this sub-expression and its closing parentheses $)^p$.
Under $P$, the set of atoms of  $\N[h^{-1}]^{(n+1)}$ is mapped to that of
$\N^{n}$, hence the last claim \cite{Lerden}.
\eproof

\noindent
The following items give moves on parentheses (and thus on planar rooted trees) to obtain all the vectors greater than a given one. 
\begin{theo}(Moves)
\label{RP}
The following holds.
\begin{enumerate}
\item {If $\vec{v}_i < \vec{w}_i $ in $\N[h^{-1}]^{n_i}$ and $\vec{v}_1 \in \N[h^{-1}]^{n_1}, \ldots, \vec{v}_m \in \N[h^{-1}]^{n_m}$, then in $\N[h^{-1}]^{n_1+ \ldots + n_m+m-1},$ 
$$\vec{v}_1 \vee \ldots \vee \vec{v}_i \vee \ldots \vee \vec{v}_m < \vec{v}_1 \vee \ldots \vee \vec{w}_i \vee \ldots \vee \vec{v}_m.$$ }
\item {If $\vec{v}:= \vec{v}_1 \vee \ldots \vee \vec{v}_m$, then for all
$\vec{w}_1 \in \N[h^{-1}]^{n_1}, \ldots, \vec{w}_p \in \N[h^{-1}]^{n_p}$,
$$ \vec{v} \vee \vec{w}_1 \vee \ldots \vee \vec{w}_p < \vec{v}_1 \vee \ldots \vee \vec{v}_m \vee \vec{w}_1 \vee \ldots \vee \vec{w}_p.$$}
\item {If $\vec{v}:= \vec{v}_1 \vee \ldots \vee \vec{v}_m$ in $\N[h^{-1}]^{n}$, then for all $0<j<m$,
$$ \vec{v} < \vec{v}_1 \vee \ldots \vee \vec{v}_{j} \vee ( \vec{v}_{j+1} \vee \ldots \vee \vec{v}_m).  $$}
\end{enumerate} 
Moreover, every vector greater than a given one can be obtained by action of Items 1, 2 and 3. 
\end{theo}
\Proof
The first claim comes from the definition of the grafting operation (see Proposition \ref{sed}) and
the partial order on $\N[h^{-1}]^n$. Consider now the inequality written in Item 2 and observe, on the left hand side, that the last coordinate of  $\vec{v}_m$ is a polynomial starting with $1+ h^{-1}+ P[h^{-1}]$, with  $P[h^{-1}]= \sum_{i \in I} h^{-i} $. Similarly, the last coordinate of  $\vec{w}_p$ is a polynomial starting with $1+ h^{-1}+ Q[h^{-1}]$,
with  $Q[h^{-1}]= \sum_{i \in J} h^{-i} $.
However, on the right hand side we have two possibilities. Either this coordinate placed in, say $i>1$, becomes $i-1 + ih^{-1}$, or $i_0 + P[h^{-1}]$, with $i_0>1$ is the label of the first coordinate of $\vec{v}_m$ in the parenthesis expression at the right hand side of this inequality. Observe the last coordinate of  $\vec{w}_p$  will not be modified.
The proof of Item 3 is complete by observing that the first coordinate $\vec{v}_{j+1}$, placed in say $k$ upgrades or remains to $k$ and that
the last coordinate of  $\vec{v}_m$, polynomial starting with $1+ h^{-1}+ Q[h^{-1}]$,
with  $Q[h^{-1}]= \sum_{i \in J} h^{-i} $ will be $1+ h^{-1}+ Q[h^{-1}] + h^{-k}$. Focus on the last claim. If $\vec{v} <  \vec{w}$, then this means the existence of say $k$ different coordinates. The case $k=1$ is 
uniquely obtained in the following case, $\ldots X_i)^p \ldots X_j)^q$,  where,
$v_i=j_0 + h^{-j_0} + \ldots $, $v_j=j_0 + h^{-j_0} + \ldots $ and $j$ being the lower index such that there is no complete expression englobing that begining in $j_0$ and endding in $i$, to $\ldots X_i)^{p-1} \ldots X_j)^q$,  --(suppression of a parenthesis)-- since in this case, $w_i=j_0' + h^{-j_0'} + \ldots $ and $w_j=j_0 + h^{-j_0} + \ldots :=v_j$. This corresponds to the move of Item 2 inside the complete expression between $j_0$ and $i$. The case $k=2$ is always obtained by Items 2 and 3. For instance,
$\ldots X_i)^p \ldots X_j)^q$, where
$v_i=i_0 + h^{-i_0} + \ldots $, $v_j=j_0 + h^{-j_0} + \ldots $ and $j_0 < i_0$ and $j$ being the lower index such that there is no complete expression englobing that begining in $i_0$ and endding in $i$, to $\ldots X_i)^{p-1} \ldots X_j)^q$, since in this case, $w_i=i_0' + h^{-i_0'} + \ldots $, with $i_0' > i_0$ and $w_j=j_0 + h^{-j_0} + \ldots + h^{i_0}$, which is the move described by Item 2 within the complete expression
starting with the most external ( in $i_0$ and closed in $j$ and Item 1 within this whole complete expression. Check the other cases to complete proof. 
\eproof
\Rk
In \cite{PalRon}, Theorem \ref{RP} was proposed as a definition to introduce a partial order on $T_n$.  In fact, this partial order comes from a natural coding of $T_n$ \textit{via} parentheses. This partial order generalises the so-called Tamari lattice on $Y_n$.
 
\noindent
It has been proved in \cite{BlassSagan} that Tamari lattices have a deep property: they are $LL$-lattices. We will now show that our generalisation of Tamari lattices are also $LL$-lattices. Let $(\mathcal{L}, <)$ be a lattice with minimal element $\hat{0}$ and maximal $\hat{1}$. In the case of a supersolvable, --\textit{i.e.}, having a maximal chain $\delta:=\hat{0}:=x_0 < x_1 < \ldots < x_{n-1} < x_n:= \hat{1}$ verifying some properties-- and semi-modular lattice, R. P. Stanley \cite{Stanley} proved that its characteristic polynomial $\chi(\mathcal{L}, x)$ factors as $\Pi(x - a_i)$, where
$a_i$ are the numbers of atoms of $\mathcal{L}$ below $x_i$ but not below $x_{i+1}$. Stanley hypotheses have been weakened by A. Blass and B. E. Sagan \cite{BlassSagan}. Let $(\mathcal{L}, <)$ be a lattice with minimal element $\hat{0}$ and maximal $\hat{1}$ equipped with a maximal chain
$\delta:=\hat{0}:=x_0 < x_1 < \ldots < x_{n-1} < x_n:= \hat{1}$. This chain induces a partition of the set of atoms of $(\mathcal{L}, <)$ into sets $A_i=\{ a \in \mathcal{A}(\mathcal{L}, <), \  a \leq x_i \ \textrm{and} \ a \nleq x_{i-1} \}$, defined for all $0 <i \leq n$ and called levels of  $\mathcal{A}(\mathcal{L}, <)$. A partial order $\lhd$ is introduced on $\mathcal{A}(\mathcal{L}, <)$, by declaring that if $a \in A_i$ and $b \in A_j$, then $a \lhd b \Leftrightarrow i < j$.
\begin{theo}[A. Blass and B. E. Sagan \cite{BlassSagan}]
If $(\mathcal{L}, <)$ is a $LL$-lattice, \textit{i.e.,} verifies:
\begin{enumerate}
\item {(Left-modularity). There exists a maximal chain $\delta:=\hat{0}:=x_0 < x_1 < \ldots < x_{n-1} < x_n:= \hat{1}$, whose all the elements are left-modular, \textit{i.e.}, for all $y,z \in (\mathcal{L}, <)$ such that $y \leq z$ and for all $i$, $y \bigvee (x_i \bigwedge z)=(y \bigvee x_i) \bigwedge z$ holds.}
\item {(Level condition). The induced partial order $\lhd$ verifies the following conditions. For all $a, b_1, \ldots b_k \in \mathcal{A}(\mathcal{L}, <)$, $a \lhd b_1\lhd \ldots \lhd b_k \Rightarrow a \nleq \bigvee_{i=1}^k b_i$,}
Then, its characteristic polynomial $\chi((\mathcal{L}, <), x)$ factors as
$\Pi_{i=1}^n (x- \card(A_i))$, for all $x \in K$.
\end{enumerate} 
\end{theo}
\begin{lemm}
Let $\vec{u} \in \N[h^{-1}]^n$. If $\vec{u}$ is left-modular, then so is $\vec{u}^\dagger$.
\end{lemm}
\Proof
Suppose $\vec{u} \in \N[h^{-1}]^n$ to be left-modular. By definition, 
for all $\vec{v}^\dagger, \vec{w}^\dagger \in \N[h^{-1}]^n$ such that $\vec{w}^\dagger \leq \vec{v}^\dagger$, the equation $\vec{w}^\dagger \bigvee (\vec{u} \bigwedge \vec{v}^\dagger)=(\vec{w}^\dagger \bigvee \vec{u}) \bigwedge \vec{v}^\dagger$ holds. By applying Proposition \ref{meetjoint} and the dendriform involution on the previous equation, we get, for all $\vec{v}, \vec{w} \in \N[h^{-1}]^n$ such that $\vec{v} \leq \vec{w}$, the equality $ (\vec{v} \bigvee \vec{u}^\dagger)\bigwedge  \vec{w}=\vec{v} \bigvee ( \vec{u}^\dagger \bigwedge \vec{w})$.
\eproof
\begin{lemm}
Fix $n>0$. Let $\vec{u} \in \N[h^{-1}]^n$. If $\vec{u}$ is left-modular, then so is $\vec{u} \vee (0) \in \N[h^{-1}]^{n+1}$.
\end{lemm}
\Proof
Fix $n>0$ and consider the map $\natural: \N[h^{-1}]^{n+1} \longrightarrow \N[h^{-1}]^{n}$, $\vec{v}:=(v_1, \ldots v_n, v_{n+1}, v_{n+2}) \mapsto \vec{v}^\natural=(v_1, \ldots v_n, 1+h^{-1} + \sum_{r \in I})$, where as usual, the set $I$ is the set of $($ remained open during this process in the expression associated with $\vec{v}^\natural$. We have $\vec{v} < \vec{w} \Rightarrow \vec{v}^\natural \leq \vec{w}^\natural$, $\bigvee(\vec{v}, \vec{w})^\natural = \bigvee(\vec{v}^\natural, \vec{w}^\natural)$ and $\bigwedge(\vec{v}, \vec{w})^\natural = \bigwedge(\vec{v}^\natural, \vec{w}^\natural)$. The first claim holds by checking all the possibilities. The second one holds
since the process to calculate the least upper bound at the coordinate $i$ depends only on the coordinates below $i$. We have just to check the last coordinate. But the last coordinate of  $\bigvee(\vec{v}^\natural, \vec{w}^\natural)$ is by definition $1+h^{-1} + \sum_{r_1 \in I_1} h^{-r_1}$, where $r_1$ is the set of all parentheses $($ open in $\bigvee(\vec{v}^\natural, \vec{w}^\natural)$ between the indices $1$ and $n$. But, by definition, the last coordinate of 
$\bigvee(\vec{v}, \vec{w})^\natural$ is $1+h^{-1} + \sum_{r \in I} h^{-r}$, where $r$ is the set of indices corresponding to all parentheses $($ open in $\bigvee(\vec{v}, \vec{w})$ between the indices $1$ and $n$. As our process to calculate the least upper bound at the coordinate $i$ depends only on the coordinates below $i$, hence $I_1=I$. The same argument holds for the last claim.

\noindent
Choose $\vec{v} \leq \vec{w}  \in \N[h^{-1}]^{n+1}$. We have to prove that 
$ (\vec{v} \bigvee \vec{u} \vee (0))\bigwedge  \vec{w}=\vec{v} \bigvee ( \vec{u} \vee (0) \bigwedge \vec{w})$ holds, which is in fact the case for the first $n$ coordinates since by hypothesis $\vec{u}$ is left-modular.
The $(n+1)^{th}$ coordinate of $\vec{u}$ is of the form $1 + h^{-1} + \sum_{r \in I} h^{-r}$, where $I$ is the set of indices corresponding to all parentheses $($ remaining open between 1 and $n$. Notice, that $u_{n+2}:=1+h^{-1}$. Let us focus on $(n+1)^{th}$ coordinate of $ (\vec{v} \bigvee \vec{u} \vee (0))\bigwedge  \vec{w}$. Suppose the $(n+1)^{th}$ coordinate of $\vec{v}$ is $i_0 + h^{-i_0} + \sum_{r \in I_v}$, where $i_0$ is the indice of the most external $)$ standing at $n+1$ and closing one $($ open in $i_0$ and $I_v$ is the set  of indices corresponding to all parentheses $($ remaining open  in the expression of $\vec{v}$ between $i_0$ and $n$. Then, by construction of the least-upper bound, the $(n+1)^{th}$ coordinate of $(\vec{v} \bigvee \vec{u} \vee (0))$ will be $i_0 + h^{-i_0} + \sum_{r \in I'}$, where $I'$ is the set  of indices corresponding to all parentheses $($ in the expression of $(\vec{v} \bigvee \vec{u} \vee (0))$ remaining open between $i_0$ and $n$. Suppose now the $(n+1)^{th}$ coordinate of $\vec{w}$ is of the form
$j_0 + h^{-j_0} + \sum_{r \in I_w}$, where $I_r$ is the set of indices corresponding to all parentheses $($ in the expression of $\vec{w}$. Then, the $(n+1)^{th}$ coordinate of $(\vec{v} \bigvee \vec{u} \vee (0))\bigwedge  \vec{w}$ will be $\min(i_0,j_0) + h^{-\min(i_0,j_0)} + \sum_{r \in I}$, where $I$ is the set  of indices corresponding to all parentheses $($ in the expression of $(\vec{v} \bigvee \vec{u} \vee (0))\bigwedge  \vec{w}$
remaining open between $\min(i_0,j_0)$ and $n$. By the same arguments,  and under the same hypotheses, we will find that the $(n+1)^{th}$ coordinate of $\vec{v} \bigvee ( \vec{u} \vee (0) \bigwedge \vec{w})$
will be $\min(i_0,j_0) + h^{-\min(i_0,j_0)} + \sum_{r \in I_0}$, where $I_0$ is the set  of indices corresponding to all parentheses $($ in the expression of $\vec{v} \bigvee ( \vec{u} \vee (0) \bigwedge \vec{w})$
remaining open between $\min(i_0,j_0)$ and $n$. But between, $\min(i_0,j_0)$ and $n$, since $\vec{u}$ is left-modular, we will get $I=I_0$, hence the equality of the $(n+1)^{th}$ coordinates under this hypotheses. Checking the other cases give the same results. The case of the $(n+2)^{th}$ coordinates are straightforward because of the equality of the coordinates between $1$ and $n+1$ and the fact that $u_{n+2}:=1+h^{-1}$ is just a particular case of what have been just explained.
\eproof
\NB
Combining these two lemmas give a lot of possibilities. If $\vec{u}$ is
left modular, then so is $(\vec{u}^{\dagger} \vee (0)) \vee (0)$ and
$ (0) \vee ( (0) \vee \vec{u})$, and so on. This will be helpful for the following theorem.
\begin{theo}
For all $n>2$, the lattice $(\N[h^{-1}]^n, <)$ is a $LL$-lattice. Moreover, its caracteristic polynomial is $\chi(n, x)=x^{(n-1)^2}(x-1)^{(n-1)}$, for all $x \in K$.
\end{theo}
\Proof
By induction, construct the maximal chain as follows. Start with $n=2$, and consider $\delta_2:=\vec{\treeAB} < \vec{\treeM} < \vec{\treeBA}$. For $n>2$, define by induction, $\delta_n := \delta_{n-1} \vee (0) < (0) \vee \vec{n-2} \vee (0) < x_1 < x_2 < \ldots < x_n:= \vec{n}$, where
for all $n$, $\vec{n}:=(1,2,3, \ldots,n, 1+h^{-1}+h^{-2}+\ldots h^{-n})$ and 
$x_1:= (0) \vee (\vec{n-2} \vee (0))$, $x_2:= (0) \vee ((0) \vee \vec{n-3} \vee (0))$, $x_3:= (0) \vee ((0) \vee (\vec{n-3} \vee (0)))$ and so on.  Using Theorem \ref{RP}, observe that if $\delta_{n-1}$ is a maximal chain, then so is $\delta_{n}$. This turns out to be the case since $\delta_2$ has three elements, 
$\delta_3$ has seven elements, are maximal (check by hand) in $(\N[h^{-1}]^2, <)$ and $(\N[h^{-1}]^3, <)$. Hence, the lenght of $\delta_{n}$ is $(n-1)^2 + (n-1) +1$. If $\delta_n$ is left-modular so is
$\delta_n \vee (0)$. We have only to prove that $(0) \vee \vec{n-2} \vee (0)$ is left-modular, the others $x_i$ being obtained by grafting of left-modular vectors, (use the Remark just above). Therefore, we have only to show that for $n>2$, $(0) \vee \vec{n-2} \vee (0):=(1,2,3,4, \ldots, n, 2+h^{-2}+ \ldots h^{-n},1+h^{-1})$ is left-modular, which does not present any difficulties.
The level condition is automatically satisfied because of the coordinate definition of the atoms. The levels $A_i$ are either empty or are singleton, as there are $n-1$ atoms, we get $(n-1)^2 + (n-1) +1 - (n-1) -1$ elements of $\delta_n$ (the minimal element does not participate to this contribution) whose levels are empty, hence the factorisation of the characteristic polynomial.
\eproof
\subsection{The free dendriform trialgebra}
We now use this partial order on the names of planar rooted trees to exhibit an associative operation. First of all, we need to describe the grafting operation on names of trees. If $\vec{v}[h^{-1}]$ is a vector on $\N[h^{-1}]^n$, $n>0$, then $\vec{v}[h^{-1}]^{\flat}$ is the vector $\vec{v}[h^{-1}]$ without its last coordinate and $v$ its lenght, \textit{i.e.,} $n+1$.
\begin{prop}
\label{sed}
Let $\vec{v}_1[h^{-1}] \in \N[h^{-1}]^{n_1}, \ldots, \vec{v}_m[h^{-1}] \in \N[h^{-1}]^{n_m}$ with $n_i \not=0$ for all $i=1, \ldots, m$. Set $\vec{v}_m[h^{-1}]:=(\vec{v}_m[h^{-1}]^{\flat}, 1+h^{-1}+Q[h^{-1}])$. Then, the name of their grafting is, 
\begin{footnotesize}
$$\vec{v}_1[h^{-1}] \vee \ldots \vee \vec{v}_m[h^{-1}]=(\vec{v}_1[h^{-1}], v_1\boxplus \vec{v}_2[h^{-1}], \ldots, (v_1+ \ldots + v_{m-1})\boxplus \vec{v}_m[h^{-1}]^{\flat}, 1+h^{-1} + (v_1+ \ldots + v_{m-1})\boxplus Q[h^{-1}]),$$
\end{footnotesize}
where $k \boxplus \vec{v}_i[h^{-1}]$ shifts its coordinates $v_{i,j}[h^{-1}]$ in the following way. If $v_{i,j}[h^{-1}]:=j$, then $k \boxplus v_{i,j}[h^{-1}]:= k+j$. If $v_{i,j}[h^{-1}]:=(j-1)+ jh^{-1} $, then $k \boxplus v_{i,j}[h^{-1}]:= (k+j-1)+ (k+j)h^{-1}$. If $v_{i,j}[h^{-1}]:=j_0+ h^{-j_0} + h^{-j_0'} +\ldots $, then $k \boxplus v_{i,j}[h^{-1}]:= (k+j_0)+ h^{-(j_0 + k)} + h^{-(j_0'+ k)} +\ldots $.
If there exists $1<i<m$, with $\vec{v}_i[h^{-1}]:=(0)$, then replace $(v_1+ \ldots + v_{i-1})\boxplus \vec{v}_i[h^{-1}]$ in the previous equality by $(v_1+ \ldots + v_{i-1}) + (v_1+ \ldots + v_{i-1}+1)h^{-1}$. In addition, $(0)\vee \vec{v}_1[h^{-1}] \vee \ldots \vee \vec{v}_m[h^{-1}]:=(1, 1 \boxplus  (\vec{v}_1[h^{-1}] \vee \ldots \vee \vec{v}_m[h^{-1}])^{\flat}, 1+h^{-1}+ (1+v_1+ \ldots + v_{m-1})\boxplus Q[h^{-1}])$, 
and $\vec{v}_1[h^{-1}] \vee \ldots \vee \vec{v}_m[h^{-1}] \vee (0):=( \vec{v}_1[h^{-1}] \vee \ldots \vee \vec{v}_m[h^{-1}], 1+h^{-1})$.
Furthermore, if $t \in T_n$ and $r \in T_m$ are trees with names $\vec{v}[h^{-1}]$ and $\vec{w}[h^{-1}]:=(1:=w_1[h^{-1}],w_2[h^{-1}], \ldots,w_{m+1}[h^{-1}])$ resp., then $t \nearrow r$, the tree with the root of $t$ glued with the most left leave of $r$, is named:
$$\vec{v}[h^{-1}] \nearrow \vec{w}[h^{-1}]= (\vec{v}, (v-1) \triangleright w_2[h^{-1}],  (v-1) \triangleright w_3[h^{-1}], \ldots, (v-1) \triangleright w_{m+1}[h^{-1}]),$$
where for all $i>2$, $(v-1) \triangleright w_i[h^{-1}]= (v-1) \boxplus w_i[h^{-1}]$, except if $w_i[h^{-1}]=1+ h^{-1} +Q[h^{-1}]$. In this case, $ (v-1) \triangleright w_i[h^{-1}]=1+ h^{-1} +  (v-1) \boxplus Q[h^{-1}]$. Similarly $t\nwarrow r$, the tree with the root of $r$ glued with the most right leave of $t$, is named:
$$\vec{v}[h^{-1}] \nwarrow \vec{w}[h^{-1}]=(\vec{v}^{\flat}, (v-1) \boxplus \vec{w}^{\flat}, v_{n+1}[h^{-1}] + (v-1) \boxplus Q[h^{-1}]),$$
if the last coordinate of $\vec{w}[h^{-1}]$ is of the form $1+ h^{-1} +  Q[h^{-1}]$. These two associative operations are extended to $\treeO$, which play the r\^ole of the unit for these two operations. Moreover, 
$ \vec{v}[h^{-1}] \nearrow \vec{w}[h^{-1}] < \vec{v}[h^{-1}] \nwarrow \vec{w}[h^{-1}],$
$(\vec{v}[h^{-1}] \nearrow \vec{w}[h^{-1}])^\dagger=\vec{w}[h^{-1}]^\dagger \nwarrow \vec{v}[h^{-1}]^\dagger $ and $(\vec{v}[h^{-1}] \nwarrow \vec{w}[h^{-1}])^\dagger=\vec{w}[h^{-1}]^\dagger \nearrow \vec{v}[h^{-1}]^\dagger $.
\end{prop}
\Proof
Use the analogy between planar trees and complete expressions on $\bra \ x_1, \ldots, x_n, x_{n+1}, ( ,)\ \ket$ to complete the proof.
\eproof
\begin{prop}
\label{tridass}
Fix $p,n,m \not=0$ and $\vec{u} \in \N[h^{-1}]^p, \ \vec{v} \in \N[h^{-1}]^n$ and $\vec{w} \in \N[h^{-1}]^m$.
Then, the binary operation $\star$ defined as follows,
$$ \vec{v} \star \vec{w} := \sum_{ \vec{v} \nearrow \vec{w} \leq \vec{u} \leq  \vec{v} \nwarrow \vec{w}} \vec{u},$$
is associative. Moreover, 
$$ \vec{u} \star \vec{v} \star \vec{w} := \sum_{  \vec{u}\nearrow \vec{v} \nearrow \vec{w} \leq \vec{t} \leq   \vec{u} \nwarrow \vec{v} \nwarrow \vec{w}} \vec{t}.$$ 
This associative operation is compatible with $(0)$, \textit{i.e.}, $\vec{v} \star (0)=\vec{v}=(0) \star \vec{v}$. Furthemore, $(\vec{v} \star \vec{w})^\dagger=\vec{w}^\dagger \star \vec{v}^\dagger$.
\end{prop}
\Proof
Keep notation of Proposition \ref{tridass}.
Observe that $\nearrow$ and $\nwarrow$ are associative operations and verify $ \vec{u} \nearrow  (\vec{v} \nwarrow  \vec{w}) = (\vec{u} \nearrow  \vec{v}) \nwarrow  \vec{w} $. 
Use Proposition \ref{sed} to check that $ (\vec{u} \nwarrow  \vec{v}) \nearrow  \vec{w} < \vec{u} \nwarrow  (\vec{v} \nearrow  \vec{w})$ and
to complete the proof.
\eproof
\begin{theo}
\label{tridalg}
The $K$-vector space $K\N[h^{-1}]^{\infty}_*:= \bigoplus_{n>0} K\N[h^{-1}]^n$ equipped with the following operations,
\begin{enumerate}
\item{$\vec{v} \succ (\vec{w}_1 \vee \ldots \vee \vec{w}_m):= \sum_{\vec{v} \nearrow  \vec{w}_1 \leq \vec{u} \leq  \vec{v} \nwarrow  \vec{w}_1}  \ \vec{u} \vee \vec{w}_2 \vee \ldots \vee \vec{w}_m$, }
\item{$(\vec{v}_1 \vee \ldots \vee \vec{v}_n) \prec \vec{w}:= \sum_{\vec{v}_n \nearrow  \vec{w} \leq \vec{u} \leq  \vec{v}_n \nwarrow  \vec{w}}  \ \vec{v}_1 \vee \ldots \vee \vec{v}_{n-1} \vee \vec{u}$, }
\item{$(\vec{v}_1 \vee \ldots \vee \vec{v}_n) \bullet (\vec{w}_1 \vee \ldots \vee \vec{w}_m):= \sum_{\vec{v}_n \nearrow  \vec{w}_1 \leq \vec{u} \leq  \vec{v}_n \nwarrow  \vec{w}_1}  \ \vec{v}_1 \vee \ldots \vee \vec{v}_{n-1} \vee \vec{u} \vee \vec{w}_2 \vee \ldots \vee \vec{w}_m$, }
\end{enumerate} 
for all $\vec{v}:=\vec{v}_1 \vee \ldots \vee \vec{v}_n$ and $\vec{w}:=\vec{w}_1 \vee \ldots \vee \vec{w}_m$ --where to ease notation, $[h^{-1}]$ has been dropped--
is a dendriform trialgebra generated by $(1,1+h^{-1}):=\vec{\treeA}$.
\end{theo}
\Proof
Fix $\vec{u},\vec{v}, \vec{w} \in K\N[h^{-1}]^{\infty}_*$.
We prove that the operation $\star$ is the sum of this three operations.
First of all observe that $\vec{v} \nwarrow \vec{w}= \vec{v}_1 \vee \ldots \vee (\vec{v}_p \nwarrow \vec{w}),$ if $\vec{v}:= \vec{v}_1 \vee \ldots \vee \vec{v}_p$. Similarly, if $\vec{w}:= \vec{w}_1 \vee \ldots \vee \vec{w}_q$, then $\vec{v} \nearrow \vec{w}= (\vec{v}  \nearrow \vec{w}_1) \vee \ldots \vee \vec{w}_q$. 

\noindent
The operation $\vec{v} \succ \vec{w}$ gives the `interval',
$$ \vec{v} \nearrow \vec{w}=( \vec{v} \nearrow \vec{w}_1) \vee \ldots \vee \vec{w}_q \leq \vec{u} \leq (\vec{v}_1 \vee \ldots \vee \vec{v}_{p-1}  \vee (\vec{v}_p \nwarrow \vec{w}_1)) \vee \vec{w}_2 \vee \ldots \vee \vec{w}_q.$$
Because of the majoration $\vec{v} \nwarrow \vec{w}$, observe that
no move corresponding to Items 1,2 or 3 of Theorem \ref{RP} can be applied on these vectors
except to the last one when only move corresponding to Item 2 can be applied given thus,
$$\vec{v}_1 \vee \ldots \vee \vec{v}_{p-1}  \vee (\vec{v}_p \nwarrow \vec{w}_1) \vee \vec{w}_2 \vee \ldots \vee \vec{w}_q. $$
This is the start of the definition of the operation $\bullet$. Indeed we get,
$$\vec{v}_1 \vee \ldots \vee \vec{v}_{p-1} \vee (\vec{v}_{p} \nearrow \vec{w}_1) \vee \vec{w}_2 \ldots \vee \vec{w}_q \leq \vec{u} \leq \vec{v}_1 \vee \ldots \vee \vec{v}_{p-1} \vee (\vec{v}_{p} \nwarrow \vec{w}_1) \vee \vec{w}_2 \ldots \vee \vec{w}_q.$$
Here again, no move can be applied to these vectors except for the last one where only the third one is authorised giving thus the start for the definition of the operation $\prec$. Indeed, Item 3
of Theorem \ref{RP} applies to $\vec{v}_1 \vee \ldots \vee \vec{v}_{p-1} \vee (\vec{v}_{p} \nearrow \vec{w}_1) \vee \vec{w}_2 \ldots \vee \vec{w}_q $
give the unique element $\vec{v}_1 \vee \ldots \vee \vec{v}_{p-1} \vee ((\vec{v}_{p} \nearrow \vec{w}_1) \vee \vec{w}_2 \ldots \vee \vec{w}_q)$.
The operation $\vec{v} \prec \vec{w}$ gives the `interval',
$$ \vec{v}_1 \vee \ldots \vee \vec{v}_{p-1} \vee ((\vec{v}_{p} \nearrow \vec{w}_1) \vee \vec{w}_2 \ldots \vee \vec{w}_q) \leq \vec{u} \leq 
\vec{v}_1 \vee \ldots \vee \vec{v}_{p-1}  \vee (\vec{v}_p \nwarrow \vec{w})=\vec{v} \nwarrow \vec{w}.$$
Therefore, the interval $[\vec{v} \nearrow \vec{w}, \vec{v} \nwarrow \vec{w}]$ can be separated into 3 disjoints subintervals.

\noindent
We have now less than seven axioms to check by using the dendriform involution. Indeed, if $(\vec{u} \prec \vec{v}) \prec \vec{w} = \vec{u} \prec (\vec{v} \star \vec{w})$ holds for all $\vec{u},\vec{v}, \vec{w} \in KN^{\infty}_*$ then so is $\vec{w}\dagger \succ (\vec{v}\dagger \succ \vec{u}\dagger) = (\vec{w}\dagger \star \vec{v}\dagger) \succ \vec{u}\dagger$. However, $(\vec{u} \prec \vec{v}) \prec \vec{w} :=(\vec{u}_1 \vee \ldots \vee \vec{u}_{m-1} \vee (\vec{u}_m \star \vec{v})) \prec \vec{w} :=\vec{u}_1 \vee \ldots \vee \vec{u}_{m-1} \vee ((\vec{u}_m \star \vec{v}) \star \vec{w} ),$ and  $\vec{u} \prec (\vec{v} \star \vec{w}):=
\vec{u}_1 \vee \ldots \vee \vec{u}_{m-1} \vee (\vec{u}_m \star  (\vec{v} \star \vec{w})).$ Apply Proposition \ref{sed} to conclude. The other axioms  also follow from Proposition \ref{sed}. The generator is $\vec{\treeA}$ since, 
$$ \vec{v}_{1} \vee \vec{v}_{2} \vee \ldots \vec{v}_{n-1} \vee \vec{v}_{n} = \vec{v}_{1} \succ \vec{\treeA} \bullet \vec{v}_{2} \bullet \vec{\treeA} \ldots  \bullet \vec{\treeA} \bullet  \vec{v}_{n-1} \vec{\treeA} \prec \vec{v}_{n},
$$
Use induction to complete the proof.
\eproof
\NB
We recover a dendriform trialgebra over rooted planar trees \textit{via}
the equivalence between trees and their names.
This three operations can be extended to $(0)$ by declaring that
$(0) \prec \vec{v}:= 0 =: \vec{v} \succ (0)$, $\vec{v} \prec (0):= \vec{v} =:  (0) \succ \vec{v}$ and $(0) \bullet \vec{v}:=0=: \vec{v}\bullet (0) $, for all
$\vec{v} \in K\N[h^{-1}]^{\infty}_*$. The expressions $(0) \prec (0)$, $(0) \bullet (0)$ and $(0) \succ (0)$ are not defined. We denote by  $K\N[h^{-1}]^{\infty}:=K\treeO \oplus \bigoplus_{n>0} K\N[h^{-1}]^n,$
the augmented dendriform trialgebra generated by $\treeA$.
\NB
The existence of a dendriform trialgebra structure over rooted planar trees has first been showed by J.-L. Loday and M. Ronco in \cite{LodayRonco}. Nevertheless, the explicit definitions of the three operations remained to be discovered. This gap was filled out in a very recent article
of P. Palacios and M. Ronco \cite{PalRon} \textit{via} permutations groups methods. Here, we have proposed a parenthesis method to obtain these explicit definitions.
From Theorem \ref{tridalg}, it is clear that any tree or name of a tree can be written in a unique way \textit{via} $\treeA$ or its name $(1,1+h^{-1})$ and $\prec, \succ, \bullet$. Such an expression for a tree $t$ is denoted by $\omega_t(\treeA)$ and is called the universal expression of $t$. For instance, $\treeM:=\treeA \bullet \treeA$ or
$\treeACA:=\treeA \succ \treeA  \prec \treeA$. To be complete, we recall the following result.
\begin{theo}[Loday-Ronco \cite{LodayRonco}]
The $K$-vector space $KT_{\infty}:=\bigoplus_{n>0} KT_n$ is the free dendriform trialgebra over the generator $\treeA$.
\end{theo}
\subsection{Arithmetree on rooted planar trees}
We extend our parenthesing presentation of the free dendriform dialgebra and its arithmetree \cite{Lerden} to the arithmetree over planar rooted trees. Though the extention is apparently more difficult to handle, it shelds light on
these arithmetrees and their associated involutive Hopf algebras as a possible natural way to extend, in a noncommutative way, our usual arithmetics over $(\Na,+,\times)$ because of Theorem \ref{extcoro}. Before, we need some preparations.
We reformulate the arithmetree \cite{Lodayarithm}
introduced by J.-L. Loday \textit{via} our lattice formulation. By a {\it{grove}}, we mean simply a non-empty subset of $T_n$ or $\N[h^{-1}]^n$, \textit{i.e.}, a disjoint union of trees with same number of leaves such that each tree appears only once. The set of groves over $T_n$ is denoted by $\TT_n$ and is of cardinal $2^{C_n}-1$. For instance in low degrees,
$$ \TT_0:=\{ \treeO \}, \TT_1:=\{ \treeA \},  \TT_2:=\{ \treeAB, \treeM,
\treeBA,  \treeAB \cup \treeBA,   \treeAB \cup \treeM, \ldots \}.$$
Similarly, we define $\Ng[h^{-1}]^n$ in the same way
and continue to call grove such a union of vectors. We set
$\TT_\infty := \{ \emptyset\} \cup \bigcup_{n \geq 0} \TT_n$ and 
$\Ng[h^{-1}]^\infty :=\{ \emptyset\} \cup \bigcup_{n \geq 0} \Ng[h^{-1}]^n$.
The idea is to convert the associative operation $\star$ in Proposition \ref{tridass} into an addition with values in groves.
\subsubsection{The dendriform addition}
\begin{defi}{[Dendriform addition \cite{Lodayarithm}]}
\label{opdef}
The \textit{dendriform addition} (associative though noncommutative) of two vectors $\vec{v}[h^{-1}]$ and $\vec{w}[h^{-1}]$ associated with some rooted planar trees is defined by:
$$ \vec{v}[h^{-1}] \ \dot{\pm} \ \vec{w}[h^{-1}] := \bigcup_{\vec{v}[h^{-1}] \nearrow \vec{w}[h^{-1}] \leq \vec{u}[h^{-1}] \leq \vec{v}[h^{-1}] \nwarrow \vec{w}[h^{-1}]}     \vec{u}[h^{-1}]. $$
This is extended to groves by distributivity of both sides, \textit{i.e.},
$ \cup_i \vec{v}_i[h^{-1}] \ \dot{\pm} \ \cup_j \vec{w}_j[h^{-1}]:=  \cup_{ij} (\vec{v}_i[h^{-1}]  \ \dot{\pm} \  \vec{w}_j[h^{-1}]),$ which has a meaning thanks to Theorem \ref{sum}.
For instance, at the level of trees: $\treeA \ \dot{\pm} \ \treeA:= \treeAB \cup \treeM \cup \treeBA$.
\end{defi}
As expected, the dendriform sum $\dot{\pm}$ splits into three operations on groves given by (the symbol $[h^{-1}]$ has been dropped to ease notation):
\begin{enumerate}
\item{$\vec{v} \vdash (\vec{w}_1 \vee \ldots \vee \vec{w}_m):= \bigcup_{\vec{v} \nearrow  \vec{w}_1 \leq \vec{u} \leq  \vec{v} \nwarrow  \vec{w}_1}  \ \vec{u} \vee \vec{w}_2 \vee \ldots \vee \vec{w}_m$, }
\item{$(\vec{v}_1 \vee \ldots \vee \vec{v}_m) \dashv \vec{w}:= \bigcup_{\vec{v}_m \nearrow  \vec{w} \leq \vec{u} \leq  \vec{v}_m \nwarrow  \vec{w}}  \ \vec{w}_1 \vee \ldots \vee \vec{w}_{m-1} \vee \vec{u}$, }
\item{$(\vec{v}_1 \vee \ldots \vee \vec{v}_n) \perp (\vec{w}_1 \vee \ldots \vee \vec{w}_m):= \bigcup_{\vec{v}_n \nearrow  \vec{w}_1 \leq \vec{u} \leq  \vec{v}_n \nwarrow  \vec{w}_1}  \ \vec{v}_1 \vee \ldots \vee \vec{v}_{n-1} \vee \vec{u} \vee \vec{w}_2 \vee \ldots \vee \vec{w}_m$, }
\end{enumerate} 
for all $\vec{v}:=\vec{v}_1 \vee \ldots \vee \vec{v}_n$ and $\vec{w}:=\vec{w}_1 \vee \ldots \vee \vec{w}_m$. For instance, $\Co_{[p]} \perp \Co_{[q]}= \Co_{[p+q]}$ for all integers $p,q \not=0$.
These operations are extended to groves by distributivity of both sides with respect to the disjoint union and verify on $\Ng[h^{-1}]^\infty$ the axioms:
\begin{footnotesize}
$$(\vec{u}[h^{-1}] \dashv \vec{v}[h^{-1}] ) \dashv \vec{w}[h^{-1}] = \vec{u}[h^{-1}] \dashv (\vec{v}[h^{-1}] \dot{\pm}  \vec{w}[h^{-1}]), \ 
(\vec{u}[h^{-1}] \vdash \vec{v}[h^{-1}] )\dashv \vec{w}[h^{-1}] = \vec{u}[h^{-1}] \vdash(\vec{v}[h^{-1}] \dashv \vec{w}[h^{-1}]), \ $$
$$(\vec{u}[h^{-1}] \dot{\pm} \vec{v}[h^{-1}] )\vdash \vec{w}[h^{-1}] = \vec{u}[h^{-1}] \vdash(\vec{v}[h^{-1}] \vdash \vec{w}[h^{-1}]), \  
(\vec{u}[h^{-1}] \vdash \vec{v}[h^{-1}] ) \perp \vec{w}[h^{-1}] = \vec{u}[h^{-1}] \vdash (\vec{v}[h^{-1}] \perp  \vec{w}[h^{-1}]), \ $$
$$ (\vec{u}[h^{-1}] \dashv \vec{v}[h^{-1}] )\perp \vec{w}[h^{-1}] = \vec{u}[h^{-1}] \perp (\vec{v}[h^{-1}] \vdash \vec{w}[h^{-1}]), \ 
(\vec{u}[h^{-1}] \perp \vec{v}[h^{-1}] )\dashv \vec{w}[h^{-1}] = \vec{u}[h^{-1}] \perp(\vec{v}[h^{-1}] \dashv \vec{w}[h^{-1}]), \  $$
$$ (\vec{u}[h^{-1}] \perp \vec{v}[h^{-1}] )\perp \vec{w}[h^{-1}] = \vec{u}[h^{-1}] \perp(\vec{v}[h^{-1}] \perp \vec{w}[h^{-1}]), $$
\end{footnotesize}
for any $\vec{u}[h^{-1}], \vec{v}[h^{-1}], \vec{w}[h^{-1}] \in \bigcup_{n>0} \Ng[h^{-1}]^n$.
The action of the unit $(0)$ is defined as follows, $(0) \vdash \vec{v}[h^{-1}]= \vec{v}[h^{-1}]=\vec{v}[h^{-1}] \dashv (0)$ and 
$(0) \dashv \vec{v}[h^{-1}] =\vec{v}[h^{-1}] \vdash (0)= \vec{v}[h^{-1}] \perp (0)= (0) \perp \vec{v}[h^{-1}] = \emptyset$ for any $\vec{v}[h^{-1}] \in \bigcup_{n>0} \Ng[h^{-1}]^n$. If $\circ \in \{ \dashv,  \vdash, \perp \}$, then $\vec{v}[h^{-1}] \circ \emptyset =  \emptyset = \emptyset \circ \vec{v}[h^{-1}]$ for any $\vec{v}[h^{-1}] \in \bigcup_{n \geq 0} \Ng[h^{-1}]^n$.
The symbols $(0) \dashv (0)$, $ (0) \vdash (0)$ and $(0) \bullet (0)$ are not defined, however $(0) \dot{\pm} (0)=(0).$
The compatibility with the dendriform involution still holds since
$(\vec{v}[h^{-1}] \dashv \vec{w}[h^{-1}])^\dagger=\vec{w}[h^{-1}]^\dagger \vdash \vec{v}[h^{-1}]^\dagger$, $(\vec{v}[h^{-1}] \vdash \vec{w}[h^{-1}])^\dagger=\vec{w}[h^{-1}]^\dagger \dashv \vec{v}[h^{-1}]^\dagger$ and $(\vec{v}[h^{-1}] \perp \vec{w}[h^{-1}])^\dagger=\vec{w}[h^{-1}]^\dagger \perp \vec{v}[h^{-1}]^\dagger$ implying that
$(\vec{v}[h^{-1}] \dot{\pm} \vec{w}[h^{-1}])^\dagger=\vec{w}[h^{-1}]^\dagger \dot{\pm} \vec{v}[h^{-1}]^\dagger$.  
The dendriform trialgebra operations are recovered \textit{via} the following trick \cite{Lodayarithm}. Set $ \N[h^{-1}]^{\infty}:= \bigcup_{n \geq 0} \N[h^{-1}]^n$.
Consider the $K$-vector space $K[\Ng[h^{-1}]^{\infty}]$ spanned by the set $\{X^{\vec{v}[h^{-1}]}, \ \vec{v}[h^{-1}] \in \N[h^{-1}]^{\infty} \}$ as a $K$-vector space
and consider the following three binary operations, $\prec, \succ, \bullet$ on groves of $K[\Ng[h^{-1}]^{\infty}]$
defined by $X^{\vec{v}[h^{-1}] \perp \vec{w}[h^{-1}]}:=X^{\vec{v}[h^{-1}]} \bullet  X^{\vec{w}[h^{-1}]}$, $X^{\vec{v}[h^{-1}] \vdash \vec{w}[h^{-1}]}:=X^{\vec{v}[h^{-1}]} \succ  X^{\vec{w}[h^{-1}]}$ and $X^{\vec{v}[h^{-1}] \dashv \vec{w}[h^{-1}]}:=X^{\vec{v}[h^{-1}]} \prec  X^{\vec{w}[h^{-1}]}$, where of course we set $X^{\cup_i \ \vec{v}_i[h^{-1}]}:= \sum_i X^{\vec{v}_i[h^{-1}]}$, $X^{\emptyset}:=0$ and $X^{(0)}:=1$. It is easy to see that $(K[\Ng[h^{-1}]^{\infty}], \prec,\succ,\bullet)$ is the free dendriform trialgebra one the generator $X^{(1,1+h^{-1})}$ augmented with the unit $1:=X^{(0)}$. But now, we have enrich our space with an arithmetics over planar rooted trees called arithmetree, like the usual polynomial algebra $K[X]:=\{X^n, \ n \in (\Na,+, \times)\}$. To complete this analogy, two things are missing. The analogue of $\times$ for usual integers and
to proof that our operations are in values in groves. For that, consider the following lemma.
\begin{lemm} 
\label{encadrtri}
Let $\vec{w}[h^{-1}] \in \N^{n+m}[h^{-1}]$. Then, there exists unique  $\vec{u}[h^{-1}] \in \N^{n}[h^{-1}]$ and $\vec{v}[h^{-1}]  \in \N^{m}[h^{-1}]$ such that: $\vec{u}[h^{-1}] \nearrow \vec{v}[h^{-1}] \leq \vec{w}[h^{-1}] \leq \vec{u}[h^{-1}] \nwarrow \vec{v}[h^{-1}].$
\end{lemm}
\Proof
Let $\vec{w}[h^{-1}] \in \N[h^{-1}]^{n+m}$. For $\vec{u}[h^{-1}]$, take the first $n$ coordinates of $\vec{w}[h^{-1}] \in \N[h^{-1}]^{n+m}$ and add an extra-coordinate whose r\^ole is to close the $($ let open in the expression associated with $\vec{u}$. This gives
a unique vector $\vec{u} \in \N[h^{-1}]^{n}$.  Consider the vector $\vec{v_1}[h^{-1}]$ defined by $\vec{v_1}[h^{-1}]:=(w_{n+1}[h^{-1}], \ldots, w_{n+m}[h^{-1}])$. Make the translation of $-(u-1)$ to obtain $\vec{v_1}[h^{-1}]-(u-1):=(w_{n+1}[h^{-1}]-u, \ldots, w_{n+m}[h^{-1}]-u)$. That is, if $w_i[h^{-1}]:=i$, then $w_i[h^{-1}]-(u-1):=i-(u-1)$.  If $w_i[h^{-1}]:=i-1 + ih^{-1}$, then $w_i[h^{-1}]-(u-1):=i-1-(u-1) + (i-(u-1))h^{-1}$ and if $w_i[h^{-1}]:=i_0 + h^{-i_0} + h^{-i_0'} + \ldots$, then $w_i[h^{-1}]-(u-1):=i_0-(u-1) + h^{-(i_0-(u-1))} + h^{-(i_0'-(u-1))} + \ldots$. In this last case, if  $i_0-(u-1) \leq 0$, then replace it by $1 + h^{-1}$ and discard any $h^{-(i_0'-(u-1))}$ with exponent less or equal than zero. This defines the vector $\vec{v} \in \N[h^{-1}]^{m}$ we are looking for.
\textit{Via} Proposition \ref{sed}, observe that:
$\vec{u}[h^{-1}] \nearrow \vec{v}[h^{-1}] \leq \vec{w}[h^{-1}] \leq \vec{u}[h^{-1}] \nwarrow \vec{v}[h^{-1}].$
\eproof

\noindent
We now simplify the proof of the following theorem.
\begin{theo}[Loday, \cite{Lodayarithm}]
 \label{sum}
The dendriform addition of two groves is still a grove:
$$  \dot{\pm}: \Ng[h^{-1}]^n \times \Ng[h^{-1}]^m \xrightarrow{} \Ng[h^{-1}]^{n+m}.$$
\end{theo}
\Proof
A priori, it is not immediate that trees appearing in the union defining the dendriform addition are all different. Nevertheless, consider the total grove $\u{n+1}:=\cup_{\vec{w}[h^{-1}] \in \N[h^{-1}]^{n+1}} \ \vec{v}[h^{-1}]$, for all $n>0$. By applying Lemma \ref{encadrtri}, observe that,
$$\u{n + 1}:=\cup_{\vec{v}[h^{-1}] \in \N[h^{-1}]^{n}} \cup_{\vec{v}[h^{-1}] \nearrow (1, 1+h^{-1})\leq \vec{w}[h^{-1}] \leq \vec{v}[h^{-1}] \nwarrow (1, 1+h^{-1})} \ \vec{w} := \u{n} \dotplus \u{1}.$$
Apply associativity of the dendrifrom addition and induction to obtain $\u{n} \ \dot{\pm} \  \u{m}:=\u{n} \ \dot{\pm} \ \u{1} \ \dot{\pm} \ \u{1} \ldots \ \dot{\pm} \ \u{1}:=\u{n+m}. $
\eproof
\subsubsection{The dendriform multiplication}
\begin{defi}{[Dendriform multiplication \cite{Lodayarithm}]}
The dendriform multiplication, denoted by $\ltimes$ of a vector $\vec{v} \in \N[h^{-1}]^n$ by
$\vec{w} \in \N[h^{-1}]^m$ consists to replace in the universal expression
of $\vec{v}$, $\omega_{\vec{v}}((1,1+h^{-1}))$, the symbols $\prec$ by $\dashv$, $\succ$ by $\vdash$ and $\perp$ by $\bullet$. The expression so obtained is still called the universal expression, is still denoted by
$\omega_{\vec{v}}((1,1+h^{-1}))$ and is in values on groves.
Hence, $\vec{v} \ltimes \vec{w}= \omega_{\vec{v}}((1,1+h^{-1})) \ltimes \vec{w} = \omega_{\vec{v}}(\vec{w})$.
Therefore, 
$  \ltimes: \Ng[h^{-1}]^n \times \Ng[h^{-1}]^m \xrightarrow{} \Ng[h^{-1}]^{nm}.$
\end{defi}
The dendriform multiplication is extended on groves by distributivity on the left with respect to the disjoint union, \textit{i.e.,} $\cup_i \vec{v}_i \ltimes \vec{v}:= \cup_{i}  \vec{v}_i \ltimes \vec{v},$ where $\vec{v}$ is a grove and $\vec{v}_i$ some planar trees. It is associative, distributive on the left with respect to the dendriform addition, $(1,1+h^{-1})$ is a unit and $(0)$ is a left neutral element by convention. For any groves $\vec{v}$ and
$\vec{w}$, $\vec{v}^\dagger  \ltimes \vec{w}^\dagger = (\vec{v} \ltimes \vec{w})^\dagger.$ 
\subsubsection{Involutive $\mathcal{P}$-Hopf algebra}
We will point out first the existence of a connected $\mathcal{P}$-Hopf algebra on $(K[\Ng[h^{-1}]^{\infty}], \prec,\succ,\bullet)$ and will show that this space can be viewed as a natural noncommutative version of $K[ \Na]:=K\{X^n, \ n \in (\Na, +, \times)\}$ equipped with the usual (commutative) arithmetics $(\Na, +, \times)$. For that, we refer to \cite{Lodayscd}. The space $(K[\Ng[h^{-1}]^{\infty}], \prec,\succ,\bullet)^{\otimes 2}$ turns out to be a dendriform trialgebra
under the following extension of the operations:
\begin{eqnarray}
(X^{\vec{a}} \otimes X^{\vec{b}}) \circ (X^{\vec{a'}} \otimes X^{\vec{b'}}) &:= &(X^{\vec{a}} \star X^{\vec{a'}}) \otimes (X^{\vec{b}} \circ X^{\vec{b'}}) \ \ \ \textrm{if} \ \ X^{\vec{b}} \otimes X^{\vec{b'}} \not= 1 \otimes 1, \\
(X^{\vec{a}} \otimes 1) \circ (X^{\vec{a'}} \otimes 1) &:= &(X^{\vec{a}} \circ X^{\vec{a'}}) \otimes 1, \ \ \ \textrm{otherwise},
\end{eqnarray}
for any $\circ \in \{\prec,\succ,\bullet\}$ and $X^{\vec{a}},X^{\vec{b}},X^{\vec{a'}},X^{\vec{b'}} \in K[\Ng[h^{-1}]^{\infty}]$. 
\begin{prop}  
There exists an involutive connected $\mathcal{P}$-Hopf algebra structure over $K[\Ng[h^{-1}]^{\infty}]$ given for all $\vec{v}_i \in  \N[h^{-1}]^{n_i}$, $1 \leq i \leq n$ by (Sweedler notation):
$$ \Delta(X^{\vec{v}_1 \vee \ldots \vee \vec{v}_m}):=X^{\vec{v}_1 \vee \ldots \vee \vec{v}_m} \otimes 1 + \sum \ X^{\vec{v}_{1, (1)}} \star \ldots \star X^{\vec{v}_{m,(1)}} \otimes X^{\vec{v}_{1, (2)}} \vee \ldots \vee X^{\vec{v}_{m, (2)}}.$$
Furthermore if $X^{\vec{v}[h^{-1}]}$ and $X^{\vec{w}[h^{-1}]}$ are primitive, then so is 
$X^{\vec{v}[h^{-1}]} \bullet X^{\vec{w}[h^{-1}]}$.
\end{prop}  
\Proof
Observe that the action of the unit $X^{(0)}:=1$ is compatible and coherent with
the axioms of a dendriform trialgebra, we obtain a connected $\mathcal{P}$-Hopf algebra on the augmented free dendriform trialgebra by applying \cite{Lodayscd}. As $(K[\Ng[h^{-1}]^{\infty}], \prec, \succ, \bullet)$ is a representation of the free dendriform trialgebra on one generator, we get $\Delta(X^{\vec{v}_1} \vee \ldots \vee X^{\vec{v}_m}):=\Delta(X^{\vec{v}_1} \succ X^{\vec{\treeA}} \bullet X^{\vec{v}_2} \bullet X^{\vec{\treeA}} \bullet \ldots \bullet X^{\vec{v}_{m-1}} \bullet X^{\vec{\treeA}} \prec X^{\vec{v}_m})=\Delta(X^{\vec{v}_1}) \succ (1 \otimes X^{\vec{\treeA}} +  X^{\vec{\treeA}} \otimes 1)\bullet  \Delta(X^{\vec{v}_2}) \bullet \ldots \bullet \Delta(X^{\vec{v}_{m-1}}) \bullet (1 \otimes X^{\vec{\treeA}} +  X^{\vec{\treeA}} \otimes 1)\prec \Delta(X^{\vec{v}_m}).$
Because the unit action vanishes on $\bullet$ and because this $\mathcal{P}$-Hopf algebra is connected, we obtain (Sweedler notation),
\begin{eqnarray*}
 X^{\vec{v}_2} \bullet X^{\vec{\treeA}} \bullet X^{\vec{v}_3} \bullet \ldots  \bullet X^{\vec{\treeA}} \bullet X^{\vec{v}_{m-1}} &:=&
X^{\vec{v}_{2,(1)}} \star \ldots \star X^{\vec{v}_{m-1,(1)}} \otimes 
X^{\vec{v}_{2,(2)}} \bullet X^{\vec{\treeA}} \bullet \ldots \bullet X^{\vec{v}_{m-1,(2)}} \\ & & +
X^{\vec{v}_2} \bullet X^{\vec{\treeA}} \bullet \ldots \bullet X^{\vec{v}_{m-1}} \otimes 1.
\end{eqnarray*}
On the other hand, the action of the unit on $\prec$ and $\succ$ implies $\Delta(X^{\vec{v}_1}) \succ (1 \otimes X^{\vec{\treeA}} +  X^{\vec{\treeA}} \otimes 1)= X^{\vec{v}_{1,(1)}} \otimes X^{\vec{v}_{1,(2)}} \succ (1 \otimes X^{\vec{\treeA}} +  X^{\vec{\treeA}} \otimes 1) = X^{\vec{v}_{1,(1)}} \otimes X^{\vec{v}_{1,(2)}} \succ X^{\vec{\treeA}} + X^{\vec{v}_1}\succ X^{\vec{\treeA}} \otimes 1$ and $  (1 \otimes X^{\vec{\treeA}} +  X^{\vec{\treeA}} \otimes 1) \prec \Delta(X^{\vec{v}_m})=(1 \otimes X^{\vec{\treeA}} +  X^{\vec{\treeA}} \otimes 1) \prec X^{\vec{v}_{m,(1)}} \otimes X^{\vec{v}_{m,(2)}}= X^{\vec{v}_{m,(1)}} \otimes X^{\vec{\treeA}} \prec X^{\vec{v}_{m,(2)}} + X^{\vec{\treeA}} \prec X^{\vec{v}_{m}} \otimes 1.$ Apply again
the action of the unit on $\bullet$ to conclude. If $K$ has an involution denoted by the bar notation, then $(K[\Ng[h^{-1}]^{\infty}], \prec,\succ,\bullet)^{\otimes 2}$ becomes an involutive dendriform trialgebra by extending the dendriform involution as follows, $(\lambda X^{\vec{v}} \otimes X^{\vec{w}})^\dagger :=\bar{\lambda} X^{\vec{v}^\dagger}  \otimes X^{\vec{w}^\dagger}.$ In this case,
the Hopf algebra over $K[\Ng[h^{-1}]^{\infty}]$  turns out to be involutive.
\eproof

\noindent
The following result shows that the free dendriform trialgebra linked to planar rooted trees can be viewed as a natural way to extend the usual arithmetics on integers.
To avoid misunderstanding between the $+$ symbol dedicated to $K$-vector spaces and the usual addition, we choose to denote it by $\perp$ and natural number between bracket. For instance, $[2] \perp [3]:= [2 \perp 3=5]$.
Consider $K\Na$, the free $K$-vector space spanned by $\Na$. 
Denote the usual addition and multiplication resp. by  $\perp, \ \times: K\Na \otimes K\Na \xrightarrow{} K\Na$. The $K$-vector space $(K\Na \otimes K\Na, \ \perp)$ becomes an associative algebra if $\perp$ is extended as follows.
\begin{enumerate}
\item {For all $m,p \not=0$ and for all $n \in \mathbb{N}$, $[n] \otimes [m] \perp [p] \otimes [0]:= 0$;}
\item {For all $n,q \not=0$ and for all $p \in \mathbb{N}$, $[n] \otimes [0] \perp [p] \otimes [q]:= 0$;}
\item {Otherwise,  $[n] \otimes [m] \perp [p] \otimes [q]:= [n \perp p] \otimes [m \perp q].$}
\end{enumerate} 
Denote by $(K\Co, \perp)$, the $K$-associative algebra generated by the corollas, \textit{i.e.,} $\Co_{[1]} := \treeA$, $\Co_{[2]} := \treeM$, $\Co_{[3]} := \treeCor$, and so on. 
\begin{theo}
\label{extcoro}
The associative algebra $(K\Na, \perp)$ has a natural structure of unital commutative and cocommutative connected Hopf algebra $\Delta_\Na: K\Na  \xrightarrow{} K\Na \otimes K\Na$ given by $\Delta_\Na([n]):= [n] \otimes [0] + [0] \otimes [n]$ for all $n \in \Na$ and $\Delta_\Na([0]):= [0] \otimes [0]$. The map $ \cdot \times [r]: (K\Na, \perp) \xrightarrow{} (K\Na, \perp)$, given by $[p] \mapsto [p \times r]$, is a Hopf algebra automorphism for all $r \in \Na$ different from zero. Moreover,  the linear map $ext: (K\Na, \perp) \xrightarrow{} (K\Co, \bullet)$, defined by $[p] \mapsto \Co_{[p]}$, for all $p \not=0, p \in  \Na$ is an isomorphism of associative algebras.
\end{theo}  
\Proof
As the $K$-vector space $(K\N, \perp)$ is a unital associative and commutative algebra with unit $[0]$,
the operation $\perp$, once extended as above, gives to $K\N \otimes K\N$ a structure of unital associative and commutative algebra with unit $[0] \otimes [0]$. The coproduct $\Delta_\Na$ turns $(K\N, \perp)$ into a connected Hopf algebra and the linear map $ \cdot \times [r]$ is an automorphism of Hopf algebras for all $r \not= 0$.
For the last claim, observe that
$\Co_{[p]} \perp \Co_{[p]} :=\Co_{[p \perp q]}$ and that $\Co_{[p]} \ltimes \Co_{[q]} :=\Co_{[p \times q]}$, for all $p,q \not=0$, where $\perp$ acting on corollas is defined in Item 3 of Definition \ref{opdef}.
\eproof

\noindent
With respect to their coproducts, any positive integers or any corollas are primitive. \textit{Via}  Theorem \ref{extcoro}, $\bullet$ plays for corollas the r\^ole of $\perp$ for natural integers.
Rooted planar trees, \textit{via} corollas, are then a possible extention of integers. The price to pay is the lost of the unit action on $\perp$. Indeed, denote by $\Co_{[0]}:= \treeO$ and augment $(K\Co, \bullet)$, by declaring that $\Co_{[0]} \perp \Co_{[p]}:=\Co_{[p]}=:\Co_{[p]} \perp \Co_{[0]}$. Extend the operation $\perp$ to $K\Co^{\otimes 2}$ like for $(K\Na, \perp)$ and keep the definition of $\ltimes$. Replace $(K\Na, \perp)$ by $(K\Co, \perp)$ in the hypotheses of  Theorem \ref{extcoro}. Then, observe the linear map $ext$ is an isomorphism of Hopf algebras.  Axioms of the dendriform trialgebra structure on planar rooted trees force the unit action to vanish on $\perp$ for the benefit of $ \star$. We keep the usual addition and multiplication structures on corollas but lose the unit.

With this point of view, the following proposition gives `numbers' which does not seem to have a classical representation. Denote by the symbol $t \nearrow_i t'$, the operation consisting to place
the tree $t$ on the $i^{th}$ leaf of the tree $t'$. 
\begin{prop}
Let $X^{\vec{v}} \in \N[h^{-1}]^n$ be a primitive element. Then, for any corolla
$\Co_{[2p-1]}$, $p>0$, and any $\lambda_1, \ldots, \lambda_{2p},$ with $\sum_{i=1}^{2p} \lambda_i=0$, the element $\sum_{i=1}^{2p} \lambda_i  X^{\vec{v} \nearrow_i \Co_{[2p-1]}}$ is a primitive element.
\end{prop}
\Proof
Observe that for any $i$, $\Delta(X^{\vec{v} \nearrow_i \Co_{[2p-1]}})=X^{\vec{v} \nearrow_i \Co_{[2p-1]}} \otimes 1 + 1 \otimes X^{\vec{v} \nearrow_i \Co_{[2p-1]}} +
X^{\vec{v} \otimes \Co_{[2p-1]}}$ holds.
\eproof

\noindent
This natural extension of integers, realised \textit{via} corollas, to rooted planar trees gives an extra motivation for developing arithmetrees (from an operadic point of view) on trees.  
In the case of planar rooted binary trees, the dendriform involution gives
\cite{Lerden} a new representation of Catalan numbers. Here,
two new representations of super Catalan numbers or Schr\"oder numbers are proposed.
\begin{prop}
Fix an integer $n>0$. Denote by $Inv(T_n):=\{ t \in T_n, \ t=t^\dagger \}$. Then, for all $n>0$, $card(Inv(T_{2n-1}))=C_n$ and for all $n \geq 0$, $card(Inv(T_{2n}))=C_n$. Moreover, setting $C_0=1$, for all $n>0$,
$$C_n = \sum_{k=0}^{n-1}(\sum_{j=1}^{n-k} \ \sum_{ i_0 + \ldots +i_{j} = n-k-j} C_{i_0} \ldots C_{i_{j-1}})C_k.$$
\end{prop}
\Proof
This is a consequence of the following description of the Super Catalan sets. Let $T_0$ be a set with one element. Recall that the sets $T_n$, $n>0$, are defined inductively by the formula,
$$ T_n:= \cup_{k=1}^{n} \bigcup_{ i_0 + \ldots +i_{k} = n-k} T_{i_0} \times \ldots \times T_{i_k}.$$ Let us show by induction that for all $n \geq 0$, $card(Inv(T_{2n}))=C_n$. It is true in small dimensions since by hand, one observes that $card(Inv(T_{0}))=1$, $card(Inv(T_{2}))=1$, $card(Inv(T_{4}))=3$ and $card(Inv(T_{6}))=11$. We fix $p:=2n$ and suppose that for all $j <n$, our claim holds. Fix $t \in Inv(T_{2n})$. It has $2n+1$ leaves and if $t:=t_1 \vee t_2 \vee \ldots \vee t_{m-1}\vee t_{m}$, then 
$t=t_{m}^\dagger \vee t_{m-1}^\dagger \vee \ldots \vee t_2^\dagger \vee t_1^\dagger$. To construct such a tree, we choose $k<n$ and fix $2k+1$ leaves. At each side of this $2k+1$ leaves, we will have $n-k$ leaves. So only the knowleges of the $n-k$ leaves and the $2k+1$ leaves are enough to construct a tree invariant by the dendriform involution, subject to the condition that the $2k+1$ leaves have to be also invariant. By induction, we know that $card(Inv(T_{2k}))=C_k$, thus $T_{2k}$ is in bijection with $T_k$, trees
with $k+1$ leaves. We have thus reduced the problem to determine the number of trees one can construct from $n-k + (k+1)=n+1$ leaves. Hence the result and the formula just above. 
Similarly, let us show by induction that for all $n \geq 0$, $card(Inv(T_{2n-1}))=C_n$. It is true in small dimensions since by hand, one observes that $card(Inv(T_{1}))=1$, $card(Inv(T_{3}))=3$, $card(Inv(T_{5}))=11$. We fix $p:=2n-1$ and suppose that for all $j <n$, our claim holds. Fix $t \in Inv(T_{2n-1})$. It has $2n$ leaves and if $t:=t_1 \vee t_2 \vee \ldots \vee t_{m-1}\vee t_{m}$, then 
$t=t_{m}^\dagger \vee t_{m-1}^\dagger \vee \ldots \vee t_2^\dagger \vee t_1^\dagger$. To construct such a tree, we choose $k<p$ and fix $2k$ leaves. At each side of this $2k$ leaves, we will have $n-k$ leaves. So only the knowleges of the $n-k$ leaves and the $2k$ leaves are enough to construct a tree invariant by the dendriform involution, subject to the condition that the $2k$ leaves have to be also invariant. By induction, we know that $card(Inv(T_{2k-1}))=C_k$, thus $T_{2k-1}$ is in bijection with $T_k$, trees
with $k+1$ leaves. We have thus reduced the problem to determine the number of trees one can construct from $n-k + (k+1)=n+1$ leaves. Hence for all $n>0$, $card(Inv(T_{2n-1}))=C_n$.
\eproof

\noindent
\textbf{Acknowledgments:}
The author would like to thank Michael Sch\"urmann, Uwe Franz, Rolf Gohm and Stefanie Zeidler for their very warm hospitality during his stay at the Institut f\"ur Mathematik und Informatik, Greifswald, Germany, where this paper has been written.

\bibliographystyle{plain}
\bibliography{These}

\end{document}